\documentclass{amsart}
\newtheorem{theorem}{Theorem}[section]
\newtheorem{lemma}[theorem]{Lemma}

\theoremstyle{definition}
\newtheorem{definition}[theorem]{Definition}

\newtheorem{proposition}{Proposition}
\theoremstyle{remark}

\numberwithin{equation}{section}

\begin {document}
\title{Group gradings on superalgebras}
\author{M.Tvalavadze}
\address{Department of Mathematics, Middlesex College, University of Western Ontario, London, ON, N6A 5B7,Canada}
\email{mtvalava@uwo.ca}
\author{T.Tvalavadze}
\address{School of Mathematics and Statistics, Carleton University, Ottawa, ON, K1S5B6, Canada}
\email{tvala@math.carleton.ca}

\date{}
\maketitle

\section{Introduction}
The problem of the classification of group gradings on algebras of
various types stems from Kac's famous paper \cite{Kac} in which he
essentially determined all group gradings of simple Lie algebras
by finite cyclic groups and described ${\Bbb Z}_k$-symmetric
spaces. A number of papers on gradings of Lie algerbas was
published by J. Patera and coauthors (see, for example,
\cite{Patera}, \cite{Hav}). Within the last decade, all group
gradings by finite abelian groups on simple finite-dimensional
associative algebras, Jordan algebras and most types of Lie
algebras were found by Y. Bahturin, M. Zaicev, I. Shestakov, and
others \cite{Sehgal, BSZ, BS2}.

In this paper we look into the structure of finite-dimensional
graded superalgebras of different types such as associative, Lie
and Jordan over an algebraically closed field of characteristic
zero.

First, we prove a few theorems about the finite-dimensional simple
associative superalgebras graded by finite abelian groups and
equipped with a superinvolution compatible with the grading.
Second, we apply these results to the classification of group
gradings on finite-dimensional simple Lie and Jordan superalgebras
of certain types.

\section{Definitions and introductory remarks}
Let $R=R_{\bar{0}} \oplus R_{\bar{1}}$ be an associative
superalgebra. A superinvolution on $R$ is a ${\Bbb Z}_2$-graded
linear map $*: R \rightarrow R$ of order 2 such that, for all
homogeneous $a,b \in R$, $(ab)^*=(-1)^{|a||b|}b^*a^*$. If $*$ is a
superinvolution on $R$, then the restriction of $*$ to $R_{\bar
0}$ is an involution on $R_{\bar 0}$. According to \cite{Rac},
$R=M_{n,m}(F)$ admits only two types of superinvolutions up to
conjugation with an
 automorphism of $R$. The first type is called  an {\it
orthosymplectic} superinvolution defined
in an appropriate basis as follows:
$$ \left(\begin{array}{cc}
             A& B\\
             C& D
             \end{array}\right)^{osp}=
   \left(\begin{array}{cc}
             I_n& 0\\
             0&Q_m
             \end{array}\right)^{-1}
    \left(\begin{array}{cc}
             A& -B\\
             C&  D
             \end{array}\right)^t
     \left(\begin{array}{cc}
             I_n& 0\\
             0& Q_m
             \end{array}\right)
$$
where $t$ denotes the ordinary transpose involution,
$Q_m=\left(\begin{array}{cc}
           0& I_m\\
        -I_m& 0
     \end{array}\right)$, and $I_r$ stands for
the identity matrix of order $r$. As a matter of convinuence, we
will write $X^{osp}=\Phi^{-1}X^{\tau}\Phi$ where
$X=\left(\begin{array}{cc}
             A& B\\
             C& D
             \end{array}\right)$,
$X^{\tau}=\left(\begin{array}{cc}
             A& -B\\
             C&  D
             \end{array}\right)^t$ and
 $\Phi=\left(\begin{array}{cc}
             I_n& 0\\
             0&Q_m
             \end{array}\right)$.
 The second type defined in the case $n=m$ is called a {\it
transpose} superinvolution and can be represented in an
appropriate basis as follows:
$$  \left(\begin{array}{cc}
            A& B\\
            C& D
            \end{array}\right)^{trp}=
            \left(\begin{array}{cc}
            D^t& -B^t\\
            C^t& A^t
            \end{array}\right).
$$

If $G$ is a finite abelian group, and $F$ is an algebraically
closed field of characteristic 0, then a $G$-grading on an algebra
$R$ is equivalent to an action of the dual group $\widehat G$ on
$R$ by automorphisms \cite{Sehgal}. Namely, any element $a$ in
$R=\oplus_{g\in G}\,R_g$ can be uniquely decomposed into the sum
of homogeneous components, $a=\sum_{g\in G} a_g$, $a_g\in R_g$.
Given $\chi\in \widehat G$ we can define
$$ \chi*a=\sum_{g\in G}\chi(g)a_g. \eqno (1)$$
It is easy to observe that the relation (1) defines a $\widehat
G$-action on $R$ by automorphisms and a subspace $V\subseteq R$ is
a graded subspace if and only if $V$ is invariant under this
action, i.e. $\widehat G*V=V$.

Let us recall the following definitions which will be used later
in this paper.

\begin{definition}
Let $R=\oplus_{g\in G}\,R_g$ be a $G$-graded superalgebra equipped
with a superinvolution $*$. If for each $g\in G$, $(R_g)^*=R_g$,
then $*$ is called {\it graded}.
\end{definition}
\begin{definition}
A grading $R=\oplus_{g\in G} R_g$ is called {\it fine} if for any
$g$ such that $R_g\ne\{0\}$, $\dim\, R_g=1$.
\end{definition}

\begin{definition}
A grading $ R=\oplus_{g\in G} R_g$ on the matrix algebra
$R=M_n(F)$ is called {\it elementary} if there exists an $n$-tuple
$(g_1,\ldots, g_n)\in G^n$ such that the matrix units $E_{ij}$,
$1\le i,j\le n$ are homogeneous and $E_{ij}\in R_g$ if and only if
$g=g^{-1}_ig_j$.
\end{definition}

This implies that  there exists a basis of matrix units
$\{E_{ij}\}$ such that each $E_{ij}$ is a homogeneous element. A
well known Cartan decomposition of a semisimple Lie algebra is an
example of elementary grading by a root system.

\section{Two main types of graded associative superalgebras}

In what follows we recall two types of $G$-graded
finite-dimensional associative superalgebras defined in \cite{BS}.
If
$\theta=(\underbrace{g_1,\ldots,g_1}_{k_1},\ldots,\underbrace{g_r,\ldots,g_r}_{k_r})$
is the ordered set of elements of $G$, then we will use the
following notation: $\theta=(g_1^{(k_1)},\ldots,g_r^{(k_r)})$.

{\bf Type $A(\theta; \bar p).$}\quad Let
$R=M_{p_1+\ldots+p_r,q_1+\ldots+q_r}(F)$ and
$\theta=(g_1^{(k_1)},\ldots,g_r^{(k_r)})$ where $g_i\ne g_j$,
$i,j=1,\ldots,r$, $k_1=p_1+q_1,\ldots, k_r=p_r+q_r$ define the
elementary $G$-grading on $R$. It follows that any $X\in R$ can be
represented in the following way:
$$ X=\left(\begin{array}{ccc}
         X_{11}  &  \ldots &   X_{1r}\\
        \vdots   &  \ddots &  \vdots \\
         X_{r1}  &  \ldots &  X_{rr}
         \end{array}\right)\eqno (2)
$$
where $X_{ij}=\left(\begin{array}{cc}
                A_{ij}& B_{ij}\\
                C_{ij}& D_{ij}
              \end{array}\right),$
$A_{ij}$ is a $p_i\times p_j$-matrix, $B_{ij}$ is a $p_i\times
q_j$-matrix, $C_{ij}$ is a $q_i\times p_j$-matrix, $D_{ij}$ is a
$q_i\times q_j$-matrix. Then $X$ is in $R_{\bar 0}$ if for any
$i,j=1,\ldots,r$
$$ X_{ij}=\left(\begin{array}{cc}
                A_{ij}&  0\\
                    0 & D_{ij}
              \end{array}\right).\eqno (3)
$$
A matrix $X$ is in $R_{\bar 1}$ if for any $i,j=1,\ldots,r$
$$ X_{ij}=\left(\begin{array}{cc}
                    0 &  B_{ij}\\
                  C_{ij}   &  0
              \end{array}\right).\eqno (4)
$$

Note that $R_{\bar 0}$ has the induced $G$-grading which is also
elementary, and $R_{\bar 0}=I_1\oplus I_2$, the sum of two
orthogonal ideals. Each of these ideals is also $G$-graded.

{\bf Type $Q(\theta; h).$}\quad Let $R=M_{n,n}(F)$ and
$$\theta=(g_1^{(k_1)},g_1h^{(k_1)},\ldots,g_{2r-1}^{(k_{2r-1})},g_{2r-1}h^{(k_{2r-1})})$$
where $h\in G$, $o(h)=2$, $k_1+k_3+\ldots+k_{2r-1}=n$  define the
elementary grading on $R$. It follows that any $X\in R$ can be
represented as follows:
$$ X=\left(\begin{array}{ccc}
         X_{11}  &  \ldots &   X_{1 2r-1}\\
        \vdots   &  \ddots &  \vdots \\
         X_{2r-1 1}  &  \ldots &  X_{2r-1 2r-1}
         \end{array}\right)\eqno (5)
$$
where $X_{ij}=\left(\begin{array}{cc}
                A_{ij}& B_{ij}\\
                C_{ij}& D_{ij}
              \end{array}\right),$
all $A_{ij}$, $B_{ij}$, $C_{ij}$ and $D_{ij}$ are matrices of size
$k_i\times k_j$. Then a matrix $X$ is in $R_{\bar 0}$ if for any
$i,j=1,3,\ldots,2r-1$
$$ X_{ij}=\left(\begin{array}{cc}
                A_{ij}&  B_{ij}\\
                B_{ij}&  A_{ij}
              \end{array}\right).\eqno (6)
$$
A matrix $X$ is in $R_{\bar 1}$ if for any $i,j=1,3,\ldots,2r-1$
$$ X_{ij}=\left(\begin{array}{cc}
                A_{ij}&  -B_{ij}\\
                B_{ij}&  -A_{ij}
              \end{array}\right).\eqno (7)
$$

In \cite{BS} the following theorem about the structure of
elementary gradings on associative superalgebras was obtained.
\begin{theorem}
Let $G$ be an arbitrary finite abelian group, $F$ an algebraically
closed field of characteristic different from 2, $R$ a $G$-graded
finite-dimensional associative superalgebra which is simple as an
associative algebra, whose gradings is elementary. Then as a
$G$-graded superalgebra $R$ is isomorphic to one of the
superalgebras $Q(\theta;h)$ or $A(\theta;\bar p)$.
\end{theorem}
\par\medskip
\noindent {\bf Example 1.}\quad Let
$R'=M_{p_1+\ldots+p_r,q_1+\ldots+q_r}(F)$ and
$$\theta'=(g_1^{(p_1)},g_2^{(p_2)},\ldots,g_r^{(p_r)},g_1^{(q_1)},g_2^{(q_2)},\ldots,g_r^{(q_r)})$$
where $g_i\ne g_j$, $i,j=1,\ldots,r$ define  the elementary
$G$-grading on $R$. This implies that any $A\in R'$ can be
represented as follows:
$$ A=\left(\begin{array}{cc}
                 X&  Y\\
                 Z&  T
              \end{array}\right),\eqno (8)
$$
where $X=(X_{ij})_{i,j=1}^r$, $X_{ij}$ is a $p_i\times p_j$
matrix; $Y=(Y_{ij})_{i,j=1}^r$, $Y_{ij}$ is a $p_i\times q_j$
matrix; $Z=(Z_{ij})_{i,j=1}^r$, $Z_{ij}$ is a $q_i\times p_j$
matrix; $T=(T_{ij})_{i,j=1}^r$, $T_{ij}$ is a $q_i\times q_j$
matrix. Then $A$ has a degree $g^{-1}_ig_j$ if the only possible
non-zero blocks of $A$ are $X_{ij}$, $Y_{ij}$, $Z_{ij}$ and
$T_{ij}$. Moreover, $A\in R_{\bar 0}$ if $Z=0$, $Y=0$ in (8), and
$A\in R_{\bar 1}$ if $X=0$, $T=0$ in (8).
\par\medskip
\noindent {\bf Example 2.}\quad Let $R'=M_{n,n}(F)$ and
$$\theta'=(g_1^{(k_1)},g_3^{(k_3)},\ldots,g_{2r-1}^{(k_{2r-1})},
g_1h^{(k_1)},g_3h^{(k_3)},\ldots,g_{2r-1}h^{(k_{2r-1})})$$ where
$g_i\ne g_j$, $i,j=1,3,\ldots,2r-1$, $h\in G$, $o(h)=2$,
$k_1+k_3+\ldots+k_{2r-1}=n$ define the elementary $G$-grading on
$R'$. Therefore, any $A\in R'$ can be represented as follows:
$$ A=\left(\begin{array}{cc}
                 X&  Y\\
                 Z&  T
              \end{array}\right),\eqno (9)
$$
where $X=(X_{ij})_{i,j=1}^r$,  $Y=(Y_{ij})_{i,j=1}^r$,
$Z=(Z_{ij})_{i,j=1}^r$, $T=(T_{ij})_{i,j=1}^r$, and all $X_{ij}$,
$Y_{ij}$, $Z_{ij}$, $T_{ij}$ are $k_i\times k_j$ matrices. Note
that $R_{\bar 0}$ is the sum of two ungraded orthogonal ideals
$I_1$ and $I_2$ where $I_1=\left\{\left(\begin{array}{cc}
                 X&  X\\
                 X&  X
              \end{array}\right)\right\},$
$I_2=\left\{\left(\begin{array}{cc}
                 Y&  -Y\\
                -Y&   Y
             \end{array}\right) \right\}$, $X$, $Y$ of order $n$.
Moreover, $A\in R_{\bar 0}$ if $X=T$, $Y=Z$ in (9), and $A\in
R_{\bar 1}$ if
 $X=-T$, $Y=-Z$ in (9).

\par\medskip

\begin{lemma}
As a $G$-graded superalgebra $R'$ defined in Example 1 is isomorphic to
a superalgebra of the type $A(\theta,\bar p)$ where
$\theta=(g_1^{(k_1)},\ldots,g_r^{(k_r)})$ where $g_i\ne g_j$,
$i,j=1,\ldots,r$, $k_i=p_i+q_i$, $\bar p=(p_1,\ldots,p_r)$.
\end{lemma}
\begin{proof}
Let $R$ be a superalgebra of the type $A(\theta,\bar p)$
represented in the canonical form (2), and $R'$ a superalgebra
defined in Example 1. In order to prove this lemma we will find  a
$G$-graded automorphism $\varphi:R\to R'$ such that
$\varphi(R_{\bar 0})=R'_{\bar 0}$, $\varphi(R_{\bar 1})=R'_{\bar
1}$. First of all, recall that the defining tuple for the
$G$-grading of $R$ has the following form:
$\theta=(g_1^{(k_1)},\ldots,g_r^{(k_r)})=(g_1^{(p_1)},g_1^{(q_1)},\ldots,
g_r^{(p_r)},g_r^{(q_r)})$. At the same time, the defining tuple
for $G$-grading of $R'$ is
$\theta'=(g_1^{(p_1)},g_2^{(p_2)},\ldots,g_r^{(p_r)},g_1^{(q_1)},g_2^{(q_2)},\ldots,g_r^{(q_r)})$.
Permutation of any two elements of $\theta$ results in permutation
of corresponding rows and columns in matrices from $R$, which is
clearly an automorphism. Therefore, applying a series of the above
permutations, $\theta$ can be brought to $\theta'$. Set $\varphi$
equal to the composition of all automorphisms that correspond to
permutations. Clearly, $\varphi$ is an automorphism.  If $A\in R$
is a matrix whose the only possible non-zero block $A_{ij}$ is  in
position $(i,j)$, and $A_{ij}=\left(\begin{array}{cc}
                 X_{ij}&  Y_{ij}\\
                 Z_{ij}&  T_{ij}
              \end{array}\right)$, then $\varphi(A)=
\left(\begin{array}{cc}
                 X&  Y\\
                 Z&  T
              \end{array}\right)$ is  of the form (8)
where the only possible non-zero blocks of $\varphi(A)$ are
$X_{ij}$, $Y_{ij}$, $Z_{ij}$ and $T_{ij}$. It follows that
$\varphi(R_g)=R'_g$. Clearly, $\varphi(R_{\bar 0})=R'_{\bar 0}$,
$\varphi(R_{\bar 1})=R'_{\bar 1}$. The lemma is proved.
\end{proof}
\begin{lemma}
As a $G$-graded superalgebra, $R'$ defined in Example 2 is
isomorphic to a superalgebra of the type $Q(\theta;h)$ where
$$\theta=(g_1^{(k_1)},g_1h^{(k_1)},\ldots,g_{2r-1}^{(k_{2r-1})},g_{2r-1}h^{(k_{2r-1})})$$
 $h\in G$, $o(h)=2$, $k_1+k_3+\ldots+k_{2r-1}=n$.
\end{lemma}
\begin{proof}
In the same way as in the previous lemma, we can show that a
series of permutations applied to $\theta$ can reduce $\theta$ to
$\theta'$. Let $\varphi$ be the composition of all automorphisms
corresponding to  permutations. As was shown above,
$\varphi(R_g)=R'_g$  and $\varphi(R_{\bar 0})=R'_{\bar 0}$,
$\varphi(R_{\bar 1})=R'_{\bar 1}$. The lemma is proved.
\end{proof}

\section{Involution fine gradings on $M_{n,m}(F)$}

In this section we investigate the case  where $R=M_{n,m}(F)$ is a
matrix superalgebra with a fine $G$-grading and superinvolution
$*$ compatible with the $G$-grading. First we recall the following
fact from \cite{BSZ}.

\begin{theorem}
Let $R=M_n(F)=\oplus_{g\in G}\, R_g$ be a matrix algebra over an
algebraically closed field of characteristic zero graded by the
group $G$ and \mbox{Supp}\,R generates G. Suppose that the
$G$-grading is fine and $*:R\mapsto R$ is a graded involution.
Then $G$ is abelian, $n=2^k$ and $R$ as a $G$-graded algebra with
involution is isomorphic to the tensor product
$R_1\otimes\ldots\otimes R_k$ where
\par\medskip
\emph{(1)} $R_1,\ldots,R_k$ are graded subalgebras stable under
the involution  $*$;

\emph{(2)} $R_1\otimes\ldots\otimes R_k$ is an
$H_1\times\ldots\times H_k$-graded algebra and any $R_i$, $1\le
i\le k$, is $H_i\cong {\Bbb Z}_2\times {\Bbb Z}_2$-graded algebra.
\end{theorem}

The simple consequence of the fact that a general grading on
$M_n(F)$ is a tensor product of a fine grading and an elementary
grading \cite{Sehgal} is the following  lemma.

\begin{lemma}
Let $R=M_n(F)$ be a $G$-graded algebra, $G$ a finite group. Then
$\dim{ R_e} \geq 1$. Moreover, $\dim{ R_e} = 1$ if and only if
$G$-grading is fine. $\Box$

\end{lemma}

\begin{theorem}
Let $R=M_{n,m}(F)$ be a non-trivial matrix superalgebra with
superinvolution $*$ over an algebraically closed field $F$ of
characteristic zero, and $G$ a finite abelian group. Then $R$
admits no fine $G$-gradings compatible with $*$.
\end{theorem}
\begin{proof}
Assume that $R$ has a fine $G$-grading compatible with the
superalgebra structure, i.e. $R=\oplus_{g\in G} R_g$, and $R_g =
(R_g \cap R_{\bar{0}}) \oplus (R_g \cap R_{\bar{1}}), $ for any $g
\in G.$ Additionally, we  assume that $R$ has a $G$-graded
superinvolution $*$, i.e.  $(R_g)^*=R_g$ for all $g\in G$. It
follows from \cite{BS} that $\dim{R_{\bar{0}}}=
\dim{R_{\bar{1}}}$, that is, $n=m$. Let us consider the even
component $R_{\bar 0}$.  In fact, $R_{\bar 0}$ has an involution
fine $G$-grading induced from $R$, i.e. $R_{\bar 0}=\oplus_{g\in
G} (R_{\bar 0})_g$, and $(R_{\bar 0})^*_g=(R_{\bar 0})_g$.

Since $R=M_{n,n}(F)$, $R_{\bar{0}} \cong I_1 \oplus I_2$ where
$I_1, I_2 \cong M_n(F)$. Hence, the following two cases may occur.

\noindent  Case 1:  Let $(I_1)^* = I_2$ and  $(I_2)^* = I_1$.
Hence $R_{\bar 0}$ is an involution simple algebra. Obviously, the
restriction of a fine $G$-grading of $R$ to $R_{\bar 0}$ is also
fine. According to \cite{BG}, any involution grading of a
non-simple involution simple algebra $R=A \oplus A^{op}$ is either
of the form
$$R_g=A_g \oplus A^{op}_g$$ or

$$R_g=\{(x,x^{\dag}) | x \in A_g  \} \oplus \{ (x,-x^{\dag}) | x \in A_{gh} \}
$$ where $A=\oplus_{g\in G} A_g$, ${}^{\dag}$ is a graded involution on $A$,  $h\in G$,
$o(h)=2$. Consequently, any non-simple involution simple algebra
admits no fine gradings since  dimension of each homogeneous
component is greater than one.  Therefore, this case is not
possible.

\noindent Case 2: Let $(I_1)^* = I_1$ and  $(I_2)^* = I_2$. Next
let $\widehat G$ be the dual group of $G$, and $\alpha:{\widehat
G}\to\mbox{Aut}\,A$ the homomorphism accompanying our grading. If
for each $\eta \in \widehat G$, $\alpha(\eta)(I_i)=I_i$, then a
fine $G$-grading of $A$ induces $G$-gradings on both ideals such
that $A_g=(I_1)_g\oplus(I_2)_g$. In particular, $A_e=(I_1)_e\oplus
(I_2)_e$, and $(I_i)_e\ne \{0\}$. This contradicts the fact that
our $G$-grading is fine. Therefore, there always exists $\xi\in
\widehat G$ such that $\alpha(\xi)(I_1)=I_2$. Hence, ${\widehat
G}=\Lambda\cup \Lambda\xi$ where $\Lambda=\{\eta\in {\widehat
G}|\alpha(\eta)(I_i)=I_i\}$ and $\xi^2\in \Lambda$. Then
$H=\Lambda^{\bot}$ where
$$ \Lambda^{\bot}=\{g\in G\,|\, \lambda(g)=1\,\, \mbox{for\, all}\,\,
\lambda\in \Lambda\}$$ \noindent is a subgroup of $G$ of order 2
and $\widehat{G/H}\cong\Lambda$. Let $H=\{e,h\}$ where $h^2=e$.
Next we can consider the induced $\overline{G}=G/H$-grading of
$A$. Let $\bar g=gH$ for any $g\in G$. Then $A_{\bar
g}=A_g+A_{gh}$. Since $\widehat{G/H} * I_i = \Lambda
* I_i =I_i$ where $i \in {1,2 }$, $I_i$ is a ${G/H}$-graded ideal.
It follows from $A_{\bar e}=(I_1)_{\bar e}\oplus (I_2)_{\bar e}$,
$(I_i)_{\bar e}\ne\{0\}$, and $\dim\,A_{\bar e}=2$ that
$\dim\,(I_i)_{\bar e}=1$. Therefore, both ${G/H}$-gradings on
$I_1$ and $I_2$ are fine.
 Hence, according to \cite{BS2}, ${G/H} \cong H_1 \times
\ldots \times H_k $ where $H_i \cong {\Bbb Z}_2 \times {\Bbb
Z}_2.$ Consequently, for each $g\in G$, either $g^2=e$ or $g^4=e$.
Therefore, $|G|=2\cdot 2^{2l}=2^{2l+1}$, for some natural number
$l$. On the other hand, according to \cite{Sehgal}, $G={\Bbb
Z}_{n_1}\times {\Bbb Z}_{n_1}\times\ldots\times {\Bbb
Z}_{n_k}\times {\Bbb Z}_{n_k}$, $n_i\in {\Bbb N}$. Moreover,
either $n_i=2$ or $n_i=4$. Therefore, $|G|=2^{2r}\cdot
4^{2s}=2^{2r+4s}$, for some natural numbers $r$ and $s$, which is
contradiction.

\end{proof}

\section{Involution elementary  gradings on $M_{n,m}(F)$}

In this section we investigate the case  where $R=M_{n,m}(F)$ is a
matrix superalgebra with an elementary $G$-grading and
superinvolution $*$ compatible with the $G$-grading.

\begin{lemma}
Let $R=M_{n,m}(F)$ be a matrix superalgebra with an elementary
$G$-grading and a superinvolution $*$ that respects this
$G$-grading. Then $R$ cannot be of type $Q(\theta, h)$.
\end{lemma}
\begin{proof}
Assume the contrary. Let $R$ be of type $Q(\theta,h)$ and $*$ a
superinvolution compatible with $G$-grading. Note that in this
case $n=m$. Therefore, $R_{\bar 0}$ is a direct sum of two
orthogonal isomorphic ungraded ideals $I_1$ and $I_2$. Since $R$
is of type $Q(\theta, h)$, both $I_1$ and $I_2$ are ungraded. As
usual,  let $\widehat G$ be the dual group for $G$, and
$\alpha:\widehat G\to \mbox{Aut}\,R$  a homomorphism accompanying
this $G$-grading. Then we can write $\widehat
G=\Lambda\cup\Lambda\xi$ where $\Lambda$ denote the set of all
automorphisms in $\widehat G$ that leave each $I_i$ stable, and
$\alpha(\xi)(I_1)=I_{2}$. Next we denote $\alpha(\xi)$ by
$\varphi$. Note that $\varphi$ commutes with $*$. The following
two cases may occur.
\par\medskip
\noindent{Case 1:} Let $*$ be  orthosymplectic. Then, $I_1$ is
involution simple under $*$ restricted to $I_1$ which is an
involution of the transpose type. In its turn, $I_2$ is also
involution simple under $*$ restricted to $I_2$ which is an
involution of the symplectic type. Clearly, both $I_1$ and $I_2$
are not isomorphic as involution simple algebras. However,
$\varphi:I_1\to I_2$ is an isomorphism of involution simple
algebras because it commutes with $*$, which is  a contradiction.
\par\medskip
\noindent{Case 2:} Let $*$ be transpose. Then, $R_{\bar 0}$ is
involution simple. Let $\varphi$ be defined by
$\left(\begin{array}{cc}
                                      0& M_2\\
                                      M_1 & 0
                                      \end{array}\right)$ where
$M_1$ and $M_2$ are non-degenerate matrices of order $n$. Since
$\varphi$ commutes with $*$, $M_1^t=(M_2)^{-1}$. Take any $X_1\in
R_{\bar 1}$ of the form $\left(\begin{array}{cc}
                                      0& B\\
                                      0 & 0
                                      \end{array}\right).$
Then $\varphi(X^*_1)=\left(\begin{array}{cc}
                                      0& 0\\
                                      -M^t_1B^tM_1 & 0
                                      \end{array}\right)$ and
$\varphi(X_1)^*=\left(\begin{array}{cc}
                                      0& 0\\
                                      M^t_1B^tM_1 & 0
                                      \end{array}\right)$. Clearly, $\varphi(X^*_1)\ne
                                      \varphi(X_1)^*$,
which is  a contradiction. Therefore, $R$ cannot be of type
$Q(\theta,h)$. The proof is complete.
\end{proof}

The following theorem deals with a superinvolution of the
orthosymplectic type compatible with a $G$-grading.

\begin{theorem} Let $R=M_{n,m}(F)$ be a non-trivial matrix superalgebra
with an elementary $G$-grading and a superinvolution $*$
compatible with  this $G$-grading. If $*$ is orthosymplectic, then
both $n,m$ are even, and $R$, as a graded superalgebra with
superinvolution, is isomorphic to $M_{n,m}(F)$ with the elementary
grading defined by an $(n+m)$-tuple
$$(g'_1,g'_2,\ldots,g'_{n+m})=(g_1^{(p_1)},\ldots,g_r^{(p_r)},g_1^{(q_1)},\ldots,g_r^{(q_r)}),$$
for $i\ne j$, $g_i\ne g_j$. Moreover, $p_1=p_2$,
$p_3=p_4$,$\ldots$, $p_{r-1}=p_r$, $p_1+\ldots+p_r=n$, $q_1=q_2$,
$q_3=q_4$, $\ldots$, $q_{r-1}=q_r$, $q_1+\ldots+q_r=m$ and the
superinvolution $X\to X^*= \Phi^{-1}X^{\tau}\Phi$ where
$$\Phi=\tiny\left(\begin{array}{cc}
                                     \Phi_0& 0\\
                                      0    & \Phi_1
                                      \end{array}\right),$$
$$\Phi_0=\tiny\emph{diag}\left\{\left(\begin{array}{cc}
                                      0& I_{p_1}\\
                                      I_{p_1}& 0
                                      \end{array}\right),\ldots,\left(\begin{array}{cc}
                                      0& I_{p_\frac{r}{2}}\\
                                      I_{p_\frac{r}{2}}& 0
                                      \end{array}\right)\right\}$$
and $$\Phi_1=\tiny\emph{diag}\left\{\left(\begin{array}{cc}
                                      0& I_{q_1}\\
                                     -I_{q_1}& 0
                                      \end{array}\right),\ldots,\left(\begin{array}{cc}
                                      0& I_{q_{\frac{r}{2}}}\\
                                      -I_{q_{\frac{r}{2}}}& 0
                                      \end{array}\right)\right\}$$
and $g_1g_2=g_3g_4=\ldots=g_{r-1}g_{r}$. Moreover, $K(R,*)$
consists of all matrices of the type
$$ \tiny \left(\begin{array}{cc}
                 A          & B\\
                \Phi^{-1}_1B^t\Phi_0& D
                \end{array}\right)
$$
where $\Phi^{-1}_0A^t\Phi_0=-A$, $\Phi^{-1}_1D^t\Phi_1=-D$ and
$H(R,*)$ consists of all matrices of the type
$$ \tiny \left(\begin{array}{cc}
                 A          & B\\
                -\Phi^{-1}_1B^t\Phi_0& D
                \end{array}\right)
$$
where $\Phi^{-1}_0A^t\Phi_0=A$, $\Phi^{-1}_1D^t\Phi_1=D$, where
$A$ is an $n\times n$-matrix, $D$ is an $m\times m$-matrix, $B$ is
an $n\times m$-matrix.
\end{theorem}
\begin{proof}
Let $R=\oplus_{g\in G} R_g$ be an elementary $G$-grading of $R$
 respected by  $*$. The following two cases may
occur.
\par\medskip
\noindent {Case 1:} If neither $n$ nor $m$ is even, then  $R$
admits no orthosymplectic superinvolutions  \cite{Rac}.
Consequently, there are no $G$-gradings  respected by  $*$.
\par\medskip
\noindent {Case 2:} Assume that $m$ is even.  Note that the
induced $G$-grading on $R_{\bar 0}$ is also elementary. Next,
consider $R_{\bar 0}=I_1\oplus I_2$ where $I_1\cong M_n(F)$ and
$I_2\cong M_m(F)$. Since $*$ is orthosymplectic, $I_1^*=I_1$ and
$I_2^*=I_2$. By Lemma 5.1, $R$ with $G$-grading is isomorphic to
$A(\theta;\bar p)$ with $\theta=(g_1^{(k_1)},\ldots,g_r^{(k_r)})$,
$g_i\ne g_j$, $k_i=p_i+q_i$, $\bar p=(p_1,\ldots,p_r)$. Applying a
graded automorphism, $R$ can be brought to the form (8). This
allows us to write $*$ in the standard form. Namely, for any
$X=\tiny\left(\begin{array}{cc}
                  A& B\\
                  C& D
                  \end{array}\right),$
                  $X^*=\Phi^{-1}X^{\tau}\Phi$ where
$X^{\tau}=\tiny\left(\begin{array}{cc}
                  A& -B\\
                  C& D
                  \end{array}\right)^t,$ and
$\Phi=\tiny\left(\begin{array}{cc}
                  \Phi_0& 0\\
                       0& \Phi_1
                  \end{array}\right),$ $\Phi_0$ symmetric,
                  $\Phi_1$ skewsymmetric.
Now we consider only the even part of  the identity component of
the $G$-grading we are dealing with:  $R_e\cap R_{\bar
0}=M_{p_1}(F)\oplus\ldots\oplus M_{p_r}(F)\oplus
M_{q_1}(F)\oplus\ldots\oplus M_{q_r}(F)$. Hence, $R_e\cap R_{\bar
0}$ is a semisimple associative algebra. Let us set $A_i$ be the
$i$-th component in the decomposition of $R_e\cap R_{\bar 0}$, and
let us write $R_e\cap R_{\bar 0}=A=A_1\oplus A_2\oplus\ldots\oplus
A_{2r}$. Next we consider the graded automorphism
$\varphi(X)=\Phi^{-1}X\Phi$. Let us show that $\varphi(A)=A$.
Every element in $A$ has the form $X=Y^t$ for some $Y\in A$. Then,
$\varphi(X)=\varphi(Y^t)=Y^*\in A$. So, $\varphi$ is indeed an
automorphism of $A$. We know that $\varphi(A_i)=A_{\sigma(i)}$ for
a suitable permutation $\sigma$. Next we can also define the inner
automorphism $\omega$ of $R$ given by the permutation matrix $S$
which permutes the blocks $A_i$ according to $\sigma$. Hence,
$\chi=\omega^{-1}\varphi$ leaves each block $A_i$ invariant,
$\chi(A_i)=A_i$. The restriction of $\chi$ to $A_i$ is an inner
automorphism of this matrix algebra. There exists a diagonal
matrix $T=\mbox{diag}(T_1,\ldots,T_{2r})$ such that
$\chi(X)=T^{-1}XT$ for $X\in A$. If $\Omega=\Phi^{-1}TS$, then
$\Omega=\mbox{diag}(\lambda_1 I_{p_1},\ldots,\lambda_{2r}I_{q_r})$
where $\lambda_i\in F$, $\lambda_i\ne 0$, $i=1,\ldots,2r$. It
follows from $\Phi=\Omega^{-1}TS$ that both
$\Phi_0=(\Phi^0_{ij})_{i,j=1,r}$ and
$\Phi_1=(\Phi^1_{ij})_{i,j=1,r}$ are block  matrices such that in
each column of blocks and in each row of blocks we have exactly
one non-zero block. Since $\Phi$ is non-degenerate, all
$\Phi^0_{ij}$ and $\Phi^1_{ij}$ must be square matrices. Changing
a basis of $R_g$ (if necessary) and permuting
$g_1^{(p_1)},g_2^{(p_2)},\ldots,g_r^{(p_r)},$
$g_1^{(q_1)},g_2^{(q_2)},\ldots,g_r^{(q_r)}$ we may assume that:
$$ \Phi_0=\mbox{diag}\left\{I_{p_{i_1}},\ldots,I_{p_{i_l}},\left(\begin{array}{cc}
                                         0& I_{p_{i_{l+2}}}\\
                                         I_{p_{i_{l+1}}}& 0
                                         \end{array}\right),
                                         \ldots,
                                         \left(\begin{array}{cc}
                                         0& I_{p_{i_r}}\\
                                         I_{p_{i_{r-1}}}& 0
                                         \end{array}\right)\right\},\eqno (10)$$
$$ \Phi_1=\mbox{diag}\left\{\left(\begin{array}{cc}
                                         0& I_{q_{j_2}}\\
                                         -I_{q_{j_1}}& 0
                                         \end{array}\right),
                                         \ldots,
                                         \left(\begin{array}{cc}
                                         0& I_{q_{j_r}}\\
                                         -I_{q_{j_{r-1}}}& 0
                                         \end{array}\right)\right\}. \eqno (11)$$
\noindent Compatibility $*$ with the $G$-grading gives us  the
following relations:
$g_{i_1}^2=\ldots=g_{i_l}^2=g_{i_{l+1}}g_{i_{l+2}}=\ldots=g_{i_{r-1}}g_{i_r}$
and $g_{j_1}g_{j_2}=\ldots=g_{j_{r-1}}g_{j_r}$ where
$(i_1,\ldots,i_r)$ and $(j_1,\ldots,j_r)$ are two permutations of
$(1,\ldots, r)$ resulting in (10) and (11), respectively. It
follows from non-degeneracy of $\Phi_0$ and $\Phi_1$ that
$p_{i_{l+1}}=p_{i_{l+2}},\ldots,p_{i_{r-1}}=p_{i_{r}}$ and
$q_{j_1}=q_{j_2},\ldots,q_{j_{r-1}}=q_{j_{r}}$. Note that $\Phi_0$
in (10) has two types of blocks. Let us call the blocks
$I_{p_{i_k}}$ the blocks of the {\bf first type} while
$\left(\begin{array}{cc}
                       0& I_{p_{i_{k}}}\\
                      I_{p_{i_{k}}}& 0
                     \end{array}\right)$ the blocks of the {\bf second
type}. We want to show that in fact  $\Phi_0$ cannot have both
blocks of the first and the second types. For clarity, we can
assume that $i_1=1$, $i_2=2$, $\ldots$, $i_r=r$, and
$g_1^2=\ldots=g_l^2=g_{l+1}g_{l+2}=\ldots=g_{r-1}g_r$. Next we let
$X_{ij}$ denote a block matrix of order $(n+m)$ whose only
possible non-zero block is in the position $(i,j)$. Consider
$X_{1,r+k}$ where $g_{j_k}=g_1$ of degree $g_1^{-1}g_{j_k}=e$
where $e$ is the identity of $G$. If we apply $*$, then
$X^*_{1,r+k}=\Phi^{-1}X_{1,r+k}^{\tau}\Phi$ will be a matrix that
has the only non-zero block in the position $(r+k-1,1)$ (or
$(r+k+1,1)$ depending on the parity of $k$), i.e.
$g^{-1}_{j_{k-1}}g_1=e$ (or $g^{-1}_{j_{k+1}}g_1=e$) , so
$g_1=g_{j_{k-1}}$ (or $g_1=g_{j_{k+1}}$), which is a contradiction
because $g_i\ne g_j$ if $i\ne j$.  In a similar way, it can be
shown that $\Phi_0$ cannot have only blocks of the first type. It
follows from non-degeneracy of $\Phi_0$ that $p_1=p_2$, $p_3=p_4$,
$\ldots$, $p_{r-1}=p_r$, i.e. $n$ is even, and
$\Phi_0=\mbox{diag}\left\{\left(\begin{array}{cc}
                                         0& I_{p_1}\\
                                         I_{p_1}& 0
                                         \end{array}\right),
                                         \ldots,
                                         \left(\begin{array}{cc}
                                         0& I_{p_{\frac{r}{2}}}\\
                                         I_{p_{\frac{r}{2}}}& 0
                                         \end{array}\right)\right\}.$

The next  purpose is to show that permuting
$g_{j_1},\ldots,g_{j_r}$ without changing the block-diagonal form
of $\Phi_1$ we can actually assume that $g_{j_1}=g_1$, $\ldots$,
$g_{j_r}=g_r$. For this, we fix any $i$, $1\le i\le r$. Then, for
some $k$, $g_{j_k}=g_i$. Consider  $X_{i,r+k}$ of degree
$g^{-1}_ig_{j_k}=e$. If we apply $*$, then
$X^*_{i,r+k}=\Phi^{-1}X_{i,r+k}^{\tau}\Phi$ will be a matrix in
the $(i-1)$-st column of blocks and $(r+k-1)$-st row of blocks,
i.e. $g^{-1}_ig_{j_k}=g^{-1}_{j_{k-1}}g_{i-1}=e$, so
$g_{i-1}=g_{j_{k-1}}$. This implies that consecutive elements of
$r$-tuple correspond to each block of $\Phi_1$ on the main
diagonal. Hence, permuting blocks of $\Phi_1$, we can assume that
$g_{j_1}=g_1$, $\ldots$, $g_{j_r}=g_r$. Likewise, it follows from
non-degeneracy of $\Phi_1$ that $q_1=q_2$, $q_3=q_4$, $\ldots$,
$q_{r-1}=q_r$, i.e. $m$ is even.

Next, applying $*$ to $X_{ij}$, we obtain the following relations:
$g_1g_2=g_3g_4=\ldots=g_{r-1}g_r$. Conversely, if the  relation
$g^2_1=g^2_2=\ldots=g^2_r$  holds true then the corresponding
$G$-grading is respected by  $*$. The proof is complete.
\end{proof}

The next theorem deals with  a superinvolution $*$ of the
transpose type compatible with the $G$-grading.

\begin{theorem} Let $R=M_{n,n}(F)$ be a non-trivial matrix superalgebra
 with an elementary $G$-grading and a
superinvolution $*$ that respects this $G$-grading. If $*$ is
 transpose, then $R$, as a
graded superalgebra with superinvolution, is isomorphic to
$M_{n,n}(F)$ with the elementary grading defined by an $2n$-tuple
$$(g'_1,g'_2,\ldots,g'_{2n})=(g_1^{(p_1)},\ldots,g_r^{(p_r)},g_{i_1}^{(q_{i_1})},\ldots,g_{i_r}^{(q_{i_r})}),$$
for $i\ne j$, $g_i\ne g_j$, $(i_1,\ldots,i_r)$ is any permutation
of $(1,\ldots,r)$, $p_1+\ldots+p_r=n$, $q_1+\ldots+q_r=n$  with
the superinvolution $X\to X^*= \Phi^{-1}X^{\tau}\Phi$ where
$$\Phi=\left(\begin{array}{cc}
               0& I_n\\
               I_n& 0
      \end{array}\right)
$$
and $g_1g_{i_1}=g_2g_{i_2}=\ldots=g_rg_{i_r}$. Moreover, $K(R,*)$
consists of all matrices of the type
$$ \left(\begin{array}{cc}
     A&  B\\
     C& -A^t
     \end{array}\right)
$$
where $A$ is any $n\times n$ matrix, $B$ is a symmetric matrix of
order $n$, $C$ is a skew-symmetric matrix of order $n$, and
$H(R,*)$ consists of all matrices of the type
$$ \left(\begin{array}{cc}
     A& B\\
     C& A^t
     \end{array}\right)
$$
where $A$ is any $n\times n$ matrix, $B$ is skew-symmetric of
order $n$, $C$ is symmetric of order $n$.
\end{theorem}

\begin{proof}

By Lemma 5.1,  $R$ is of the type $A(\theta;\bar p)$. Arguing in
the same way as in the proof of Theorem 5.2, $R$ can be reduced to
(8) while $\theta$ takes the form
$$(g_1^{(p_1)},\ldots,g_r^{(p_r)},g_1^{(q_1)},\ldots,
g_r^{(q_r)}).$$ This allows us to represent $*$ as follows: for
any $X=\tiny\left(\begin{array}{cc}
                  A& B\\
                  C& D
                  \end{array}\right),$
                  $X^*=\Phi^{-1}X^{\tau}\Phi$ where
$X^{\tau}=\tiny\left(\begin{array}{cc}
                  A& -B\\
                  C& D
                  \end{array}\right)^t$ and
$\Phi=\tiny\left(\begin{array}{cc}
                  0& \Phi_1\\
                  \Phi_0& 0
                 \end{array}\right),$ $\Phi_0^t=\Phi_1$.
Notice that $R_e\cap R_{\bar 0}=M_{p_1}(F)\oplus\ldots\oplus
M_{p_r}(F)\oplus M_{q_1}(F)\oplus\ldots\oplus M_{q_r}(F)$. Let us
set $A_i$ be the $i$-th component in decomposition of $R_e\cap
R_{\bar 0}$, and let us write $R_e\cap R_{\bar 0}=A=A_1\oplus
A_2\oplus\ldots\oplus A_{2r}$. Using the graded automorphism
$\varphi(X)=\Phi^{-1}X\Phi$ that permutes components of $R_e\cap
R_{\bar 0}$ it can be shown that $\Phi$ can be brought to

$$\Phi=\left(\begin{array}{cc}
               0& I_n\\
               I_n& 0
      \end{array}\right)
$$
while
$\theta=(g_1^{(p_1)},\ldots,g_r^{(p_r)},g_{i_1}^{(q_{i_1})},\ldots,g_{i_r}^{(q_{i_r})})$.
Let $X_{kl}$ stand for a matrix of order $2n$  whose only possible
non-zero block is in the position $(k,l)$. Next we are able to
compute $X_{kl}^*$ where $1\le k \le r$, $1\le l\le r$. It appears
that $X^*_{kl}=X_{i_l i_k}$. It follows that
$g_k^{-1}g_l=g^{-1}_{i_l}g_{i_k}$, $g_{i_l}g_l=g_{i_k}g_k$. Then,
$X^*_{k i_l}=X_{l i_k}$ where $1\le k \le r$, $1\le l\le r$. It
follows that $g^{-1}_k g_{i_l}=g^{-1}_l g_{i_k}$. Hence
$g_{i_l}g_l=g_{i_k}g_k$. Therefore we obtained  the following
relations $g_1g_{i_1}=g_2g_{i_2}=\ldots=g_rg_{i_r}$. Conversely,
if the relations  $g_1g_{i_1}=g_2g_{i_2}=\ldots=g_rg_{i_r}$ hold
true, then the corresponding $G$-grading is respected by $*$. The
proof is complete.
\end{proof}

\section{General involution gradings on $M_{n,m}(F)$}

Now we are ready to describe all group gradings by finite abelian
groups on $M_{n,m}(F)$ stable under a superinvolution $*$. Like in
the case of graded associative algebras, a general grading appears
to be  a `mixture' of an elementary grading and  a fine grading.

\begin{definition}
Let $A$ be a superalgebra with a superinvolution $*$. Then $A$ is
said to be involution simple (or $*$-simple) if $A$ has no
non-trivial ideals stable under $*$.
\end{definition}

First we recall the following fact from \cite{BTT}.
\begin{proposition}
{\it Any finite-dimensional involution simple superalgebra with
superinvolution $*$ over an algebraically closed field of
characteristic different from 2 is isomorphic to one of the
following:
\begin{enumerate}
 \item $R=M_{n,m}(F)$ with the orthosymplectic or transpose
involution.

    \item  $R=M_{n,m}(F)\oplus M_{n,m}(F)^*$ with the ordinary
exchange involution.

    \item  $R= Q(n)\oplus Q(n)^* $ with the ordinary exchange
involution.
\end{enumerate}}
\end{proposition}

\begin{lemma}
Let $B=B_{\bar 0}+B_{\bar 1}$ be a non-simple $*$-simple
superalgebra where $*$ is a superinvolution of $B$. Then there is
an idempotent $f\in H(B,*)$ such that $f=f_0+f_1$ where $f_0\in
B_{\bar 0}$, $f_1\in B_{\bar 1}$ and $f_1\ne 0$.
\end{lemma}
\begin{proof}
According to the above classification, $B$  is isomorphic to
either $M_{n,m}(F)\oplus M_{n,m}(F)^*$ or $Q(n)\oplus Q(n)^*$.
Next we consider each of these two cases separately.

Case 1: Let $B=A\oplus A^*$ where $A=M_{n,m}(F)$ with the exchange
superinvolution. Then $f\in H(B,*)$ if and only if $f=(a,a)$.
Moreover, $f$ is an idempotent of $B$ if and only if $a^2=a$. If
we choose any non-trivial idempotent of $A$ such that $a=a_0+a_1$
and $a_1\ne 0$, then $f=(a,a)$ satisfies  all the required
conditions. For instance, we can set $a=E_{1,n+1}+E_{n+1,n+1}$.

Case 2:  Let $B=A\oplus A^*$ where $A=Q(n)$ with the exchange
superinvolution. Note that
$$ A=\left\{ \left[\begin{array}{cc}
                X& Y\\
                Y& X
                \end{array}\right]
    \right\},
$$
\noindent where $X$, $Y$ are matrices of order $n$. Consider
$a=\left[\begin{array}{cc}
                \frac{1}{2}I_n& \frac{1}{2}I_n \\
                 \frac{1}{2}I_n&  \frac{1}{2}I_n
                \end{array}\right]\in A$ where $I_n$ denotes the identity matrix of order $n$. Then $a^2=a$, $a=a_0+a_1$ and
$a_1\ne 0$. It is easily seen that $f$ is required. The lemma is
proved.
\end{proof}

\begin{lemma}
If $R=R_{\bar 0}+R_{\bar 1}$ is a superalgebra with an
antiautomorphism $\varphi$, then for any $x,y\in R$,
$$\varphi(xRy)\subseteq \varphi(y)R\varphi(x).\eqno (12)$$
\end{lemma}
\begin{proof}
Direct computations show that relation (12) holds for any
homogeneous $x$ and $y$. Therefore, by linearity of $\varphi$, it
also holds for arbitrary $x$ and $y$.

\end{proof}

As a simple consequence of this lemma we have

\begin{lemma}
Let $R$ be a superalgebra with an antiautomorphism $\varphi$, and
$f\in R$  such that $\varphi(f)=f$. Then $R'=fRf$ is a
subsuperalgebra of $R$ stable under $\varphi$.
\end{lemma}

The proof of the following lemma can be found in \cite{BS3}

\begin{lemma}
Let $R=M_n=\oplus_{g\in G} R_g$ be a matrix algebra with an
elementary $G$-grading. If $R_e=A_1\oplus A_2$ is the sum of two
simple components, then there exists $g\in G$, $g\ne e$, such that
$A_1RA_2\subseteq R_g$. $\Box$
\end{lemma}

\begin{lemma}

Let $R=C\otimes D=\oplus_{g\in G} R_g$ be a $G$-graded matrix
superalgebra with an elementary grading on $C$, and a fine grading
on $D$ over an algebraically closed field $F$ of characteristic
not 2. Let $\varphi: R\to R$ be an antiautomorphism on $R$
preserving $G$-grading and $\sigma:R\to R$ be an automorphism of
order 2 of $R$ defining a superalgebra structure on $R$. Let also
$\varphi$ act as a superinvolution on $R_e$. Then
\par\medskip
\emph{(1)} $C_e\otimes I$ is $\varphi$-stable and $\sigma$-stable
where $I$ is the unit of $D$ and hence $\sigma$ induces a
$\mathbb{Z}_2$-grading on $C_e$ and $\varphi$ induces a
superinvolution $*$ on $C_e$ compatible with
$\mathbb{Z}_2$-grading.
\par\medskip
\emph{(2)} there are $*$-subsuperalgebras $B_1,\ldots,
B_k\subseteq C_e$ such that $C_e=B_1\oplus\ldots\oplus B_k$, and
$B_1\otimes I,\ldots,B_k\otimes I$ are $\varphi$-stable and
$\sigma$-stable.
\par\medskip
\emph{(3)} if $e_i$ is the identity of $B_i$, then $D_i=e_i\otimes
D$ is $\varphi$-stable and $\sigma$-stable.
\par\medskip
\emph{(4)} $e_iCe_i\otimes I$ is $\varphi$- and $\sigma$-stable
for each $e_i$
\par\medskip
 \emph{(5)} the centralizer of $R_e=C_e\otimes I$ in
$R$ can be decomposed as $Z_1D_1\oplus\ldots\oplus Z_kD_k$ where
$Z_i=Z'_i\otimes I$, $Z'_i$ is the center of $B_i$.
\end{lemma}

\begin{proof}

It follows from \cite{Sehgal} that the identity component $R_e$
equals to $C_e\otimes I$. Since $R_e$  is $\varphi$- and
$\sigma$-stable, both $\varphi$ and $\sigma$ induce a
superinvolution $*$  and a superalgebra structure on $C_e$. Both
structures are compatible with each other.

Since $C_e$ is semisimple, it is the direct sum of simple
subalgebras,
$$ C_e=A_1\oplus\ldots\oplus A_l. $$
If for some $i$, $1\le i\le l$, $\sigma(A_i)=A_j$ where $i\ne j$,
then it is easily seen that $A'_i=A_i+A_j$ is $\sigma$-stable.
Therefore, $C_e$ can be written as a direct sum of $\sigma$-stable
superalgebras
$$ C_e=A'_1\oplus\ldots\oplus A'_s.$$
Next, if for some $i$, $1\le i\le s$, $(A'_i)^*=A'_j$ where $i\ne
j$, then $B_i=A'_i+A'_j$ is $*$-stable. Finally, $C_e$ can be
written as a direct sum of $*$-simple superalgebras
$$ C_e=B_1\oplus\ldots\oplus B_k. $$
Therefore, (1) and (2) are proved.

Next we fix $1\le i\le k$, and consider $R'=(e_i\otimes
I)(C\otimes D)(e_i\otimes I)=e_iCe_i\otimes D$ where $e_i$ is the
identity of $B_i$. Since $\varphi(e_i\otimes I)=e_i\otimes I$, by
Lemma 6.4,  $R'$ is $\varphi$-stable. Due to $\sigma(e_i\otimes
I)=e_i\otimes I$, $R'$ is also $\sigma$-stable.

 To prove (3) we
consider the following three cases.

\noindent {Case 1:} Let $B_i$ be of the type $M_{r,s}(F)$. Then
$$e_iCe_i=B_i, \eqno (13)$$
and $e_iCe_i\otimes I=B_i\otimes I$. Hence, $e_iCe_i\otimes I$ is
$\varphi$- and $\sigma$-stable. Since $e_i\otimes D$ is a
centralizer of $e_iCe_i\otimes I$ in $R'$, it is also $\varphi$-
and $\sigma$-stable.

\noindent {Case 2:} Let $B_i=A\oplus A^{*}$ where $A=M_{r,s}(F)$.
Denote the identity of $A$ by $\varepsilon_i $. Then,
$\varepsilon^*_i$ is the identity of $A^{*}$, and
$e_i=\varepsilon_i+\varepsilon^*_i$. Note that $$e_iCe_i\otimes
I=\varepsilon_iC\varepsilon_i\otimes
I+\varepsilon_iC\varepsilon^*_i\otimes
I+\varepsilon^*_iC\varepsilon_i\otimes
I+\varepsilon^*_iC\varepsilon^*_i\otimes I. \eqno (14)$$ Next we
want to prove that both $\varphi$ and $\sigma$ permute the terms
of (14) leaving $e_iCe_i\otimes I$ invariant. Without any loss of
generality we consider just one term of the form
$\varepsilon_iC\varepsilon^*_i\otimes I$. Since
$\varepsilon_iC\varepsilon^*_i\otimes I=(\varepsilon_i\otimes
I)(C\otimes I)(\varepsilon^*\otimes I)$, by (12),
$\varphi(\varepsilon_iC\varepsilon^*_i\otimes
I)\subseteq(\varepsilon_i\otimes I)(C\otimes
D)(\varepsilon^*\otimes I)=\varepsilon_iC\varepsilon^*_i\otimes D$
and $\sigma(\varepsilon_iC\varepsilon^*_i\otimes
I)\subseteq(\varepsilon_i\otimes I)(C\otimes
D)(\varepsilon^*\otimes I)=\varepsilon_iC\varepsilon^*_i\otimes
D$.

By Lemma 6.5, there exists $g\in G$, $g\ne e$ such that
$\varepsilon_i C\varepsilon^*_i\subseteq C_g$. Hence,
$\varepsilon_i C\varepsilon^*_i\otimes I\subseteq  R_g$.
Consequently, $\varphi(\varepsilon_i C\varepsilon^*_i\otimes
I)\subseteq  R_g$ and $\sigma(\varepsilon_i
C\varepsilon^*_i\otimes I)\subseteq  R_g.$ Next we take a
homogeneous $x\in \varepsilon_i C\varepsilon^*_i$ of degree $g$
and a homogeneous $y\in D$ of degree $h$ such that $x\otimes y \in
R_g$. Then $\mbox{deg}\, (x\otimes y)=gh=g$, $h=e$. This implies
$y=\lambda I$, $\lambda\in F$ for any $x\otimes y\in R_g\cap
\varepsilon_i C\varepsilon^*_i\otimes D$.  It follows that $
R_g\cap \varepsilon_i C\varepsilon^*_i\otimes D\subseteq
\varepsilon_i C\varepsilon^*_i\otimes I$.

As a result,
 $\varphi(e_iCe_i\otimes I)=e_iCe_i\otimes I $ and $\sigma(e_iCe_i\otimes I)=e_iCe_i\otimes
 I$, that is, $e_iCe_i\otimes I$ is
$\varphi$- and $\sigma$-stable. It follows that $e_i\otimes D$,
the centralizer of $e_iCe_i\otimes I$ in $R'$, is also $\varphi$-
and $\sigma$-stable.

\noindent {Case 3:} Let  $B_i=Q(s_i)\oplus Q(s_i)^{*}$. Since
$Q(s_i)=I_1\oplus I_2$ where $I_1$, $I_2$ are simple ideals
isomorphic to  $M_{s_i}(F)$, $B_i=(I_1\oplus I_2)\oplus
(I^*_1\oplus I^*_2)$. Let $\varepsilon_i$, ${\widehat
\varepsilon}_i$, $\varepsilon^*_i$, $\widehat \varepsilon^*_i$  be
the identities of $I_1,$ $I_2$, $I^*_1$, $I^*_2$, respectively.
Then we note that $\sigma(\varepsilon_i\otimes I)={\widehat
\varepsilon}_i\otimes I$. We have that
$$e_i=\varepsilon_i+\varepsilon^*_i+{\widehat \varepsilon}_i+{\widehat
\varepsilon}^*_i. \eqno (15)$$

Therefore, $e_iCe_i\otimes I=N_1\otimes I +N_2\otimes I+N_3\otimes
I+N_4\otimes I$ where
$N_1=\varepsilon_iC\varepsilon_i+\varepsilon^*_iC\varepsilon_i+\varepsilon_iC\varepsilon^*_i+\varepsilon^*_iC\varepsilon^*_i$,
$N_2=\varepsilon_iC\widehat\varepsilon_i+\varepsilon^*_iC\widehat\varepsilon_i+\varepsilon_iC\widehat\varepsilon^*_i+\varepsilon^*_iC\widehat\varepsilon^*_i$,
$N_3=\widehat\varepsilon_iC\varepsilon_i+\widehat\varepsilon^*_iC\varepsilon_i+\widehat\varepsilon_iC\varepsilon^*_i+\widehat\varepsilon^*_iC\varepsilon^*_i$,
and
$N_4=\widehat\varepsilon_iC\widehat\varepsilon_i+\widehat\varepsilon^*_iC\widehat\varepsilon_i+\widehat\varepsilon_iC\widehat\varepsilon^*_i+\widehat\varepsilon^*_iC\widehat\varepsilon^*_i$.

Arguing in the same way as in the second case, we can prove that
$\varphi(N_1\otimes I)=N_1\otimes I$ and $\varphi(N_4\otimes
I)=N_4\otimes I$. Now we consider $N_2$ and $N_3$. Suppose that
the elementary grading on $B_iCB_i$ induced from $C$ is defined by
$(g_1,g_2,g_3,g_4)$. It is easy to see that $\deg
(\varepsilon_iC\widehat\varepsilon_i)=g^{-1}_1g_3,$ $\deg
(\varepsilon_iC\widehat\varepsilon^*_i)=g^{-1}_1g_4,$ $\deg
(\varepsilon^*_iC\widehat\varepsilon_i)=g^{-1}_2g_3,$ $\deg
(\varepsilon^*_iC\widehat\varepsilon^*_i)=g^{-1}_2g_4,$ $\deg
(\widehat\varepsilon_iC\varepsilon_i)=g^{-1}_3g_1,$ $\deg
(\widehat\varepsilon^*_iC\varepsilon_i)=g^{-1}_4g_1,$ $\deg
(\widehat\varepsilon^*_iC\varepsilon^*_i)=g^{-1}_4g_2,$ $\deg
(\widehat\varepsilon_iC\varepsilon^*_i)=g^{-1}_3g_2.$

Let us take, for example, the first term
$\varepsilon_iC\widehat\varepsilon_i$ of $N_2$. Then
$\varphi(\varepsilon_iC\widehat\varepsilon_i\otimes I)\subseteq
R_{g^{-1}_1g_3}$. On the other hand, by (12),
$\varphi(\varepsilon_iC\widehat\varepsilon_i\otimes I)\subseteq
\widehat\varepsilon^*_iC\varepsilon^*_i\otimes D\subseteq
C_{g^{-1}_4g_2}\otimes D$. If we take $x\in C_{g^{-1}_4g_2}$ and a
homogeneous $y\in D$ of degree $h$ such that $\deg (x\otimes
y)=g^{-1}_1g_3$, then $g^{-1}_1g_3=g^{-1}_4g_2h$, that is,
$h=g^{-1}_1g_4g^{-1}_2g_3$. This implies that
$$\varphi(\varepsilon_iC\widehat\varepsilon_i\otimes I)\subseteq
\widehat\varepsilon^*_iC\varepsilon^*_i\otimes D_h. \eqno(16)$$
Similarly, we can show that for each term of $N_2$ there should be
$D_h$ on the right-hand side of (16). Hence $\varphi(N_2\otimes
I)=N_3\otimes D_{h}$.

Next we take the first term $\widehat\varepsilon_iC\varepsilon_i$
of $N_3$. Likewise we can show that $\varphi(\widehat
\varepsilon_iC\varepsilon_i\otimes I)\subseteq
\varepsilon^*_iC\widehat\varepsilon^*_i\otimes D_{h^{-1}}$. Hence
for each term of $N_3$ we have $D_{h^{-1}}$ on the right hand
side, and $\varphi(N_3\otimes I)=N_2\otimes D_{h^{-1}}$.

Notice that the centralizer of $(N_2+N_3)\otimes I$ in $R'$ is
$e_i\otimes D$. Hence, the centralizer of
$\varphi((N_2+N_3)\otimes I)=N_3\otimes D_{h^{-1}}+N_2\otimes D_h$
in $R'$ is $\varphi(e_i\otimes D)$. Next we take any $x\otimes y$
in $\varphi(e_i\otimes D)$. Since $x\otimes y$ lies in the
centralizer of $N_3\otimes D_{h^{-1}}+N_2\otimes D_h$, $x\otimes
y$ commutes with each element in $N_3\otimes D_{h^{-1}}$ and
$N_2\otimes D_h$. Therefore, $x$ commutes with any matrix in $N_3$
and $N_2$. Direct computations show that $x=e_i$, and
$\varphi(e_i\otimes D)=e_i\otimes K$ where $K$ is a subspace of
$D$. By dimension arguments, $K=D$, and $\varphi(e_i\otimes D)=
e_i\otimes D$.

To prove that $e_i\otimes D$ is $\sigma$-stable, we represent
$e_iCe_i\otimes I$ as follows: $e_iCe_i\otimes I=N'_1\otimes I
+N'_2\otimes I+N'_3\otimes I+N'_4\otimes I$ where
$N'_1=\varepsilon_iC\varepsilon_i+\widehat\varepsilon_iC\varepsilon_i+\varepsilon_iC\widehat\varepsilon_i+\widehat\varepsilon_iC\widehat\varepsilon_i$,
$N'_2=\varepsilon_iC\varepsilon^*_i+\widehat\varepsilon_iC\varepsilon^*_i+\varepsilon_iC\widehat\varepsilon^*_i+\widehat\varepsilon_iC\widehat\varepsilon^*_i$,
$N'_3=\varepsilon^*_iC\varepsilon_i+\widehat\varepsilon^*_iC\varepsilon_i+\varepsilon^*_iC\widehat\varepsilon_i+\widehat\varepsilon^*_iC\widehat\varepsilon_i$,
and
$N'_4=\varepsilon^*_iC\varepsilon^*_i+\widehat\varepsilon^*_iC\varepsilon^*_i+\varepsilon^*_iC\widehat\varepsilon^*_i+\widehat\varepsilon^*_iC\widehat\varepsilon^*_i$.
In the same way as above, we can show that $e_i\otimes I$ is
$\sigma$-stable. Hence (3) is proved.

To prove (4) we note that the centralizer $Z_{R'}(e_i\otimes D)$
of $e_i\otimes D$ in $R'$ is equal to $e_iCe_i\otimes I$. Since
$e_i\otimes D$ is $\varphi$- and $\sigma$-stable, so is its
centralizer.

To prove (5) we note that the centralizer $Z$ of $C_e$ in $C$ is
equal to $Z'_1\oplus\ldots\oplus Z'_k$ where $Z'_i$ is the center
of $B_i$ and the centralizer of $R_e$ in $R$ coincides with
$Z\otimes D=Z_1D_1\oplus\ldots\oplus Z_kD_k$ where
$Z_i=Z'_i\otimes I$ and $D_i=e_i\otimes D$.

Our proof is complete.

\end{proof}

In the following lemma we look into the structure of a full matrix
algebra graded by a group $G$ and equipped with an automorphism of
order 2 compatible with the $G$-grading.

\begin{lemma}
Let $R=M_n(F)=C\otimes D=\oplus_{g\in G} R_g$ be a $G$-graded
matrix algebra with an elementary grading on $C$, and a fine
grading on $D$   over an algebraically closed field $F$ of zero
characteristic. Let $\sigma$ be an automorphism of order 2
compatible with the $G$-grading. Then there are subalgebras
$B_1,\ldots, B_k\subseteq C_e$ such that
$C_e=B_1\oplus\ldots\oplus B_k$, and $B_1\otimes
I,\ldots,B_k\otimes I$ are $\sigma$-stable. Moreover, either all
$B_i$, $1\le i\le k$, are non-simple with respect to $\sigma$ or
both $C\otimes I$ and $I\otimes D$ are $\sigma$-stable.
\end{lemma}
\begin{proof}
Again we consider the decomposition $R=C\otimes D$ where $C$
carries an elementary grading and $D$ carries a fine grading. As
before, we decompose $R_e=C_e\otimes I$ as the sum of
$\sigma$-invariant simple subalgebras $B_1,\ldots,B_k$. Each such
subalgebra is isomorphic to either a matrix algebra $B_i\cong
M_{r_i}$ or  the sum of two simple subalgebras: $B_i\cong
M_{r_i}\oplus M_{r_i}$, interchangeable by $\sigma$. If we
restrict ourselves to $R_e$ then, after  a graded conjugation, we
may assume that $\sigma|_{R_e}$ is the conjugation by a matrix
$S=S_1+\ldots+S_k$ where each $S_i\in B_i$ is either of the form
$S_i=\mbox{diag}\,\{I_{p_i},-I_{q_i}\}$ where $p_i+q_i=r_i$ if
$B_i$ is simple or is of the form $\left(\begin{array}{cc}
        0& I_{r_i}\\
       I_{r_i}& 0
       \end{array}\right)$, if $B_i$ is not simple. Let $\sigma$ be given by
$\sigma(X)=\Omega^{-1}X\Omega$.  Now if we consider conjugation by
$\Omega S^{-1}$ then this will be trivial on $R_e$, that is,
$\Omega S^{-1}$ is an element of the centralizer $Z$ of $R_e$ in
$R$. Let $e_i$ be the identity element of $B_i$, $Z_i$ the center
of $B_i$. Then $Z=Z_1\otimes D\oplus\ldots\oplus Z_k\otimes D$.
This allows us to write $\Omega$ in the following way:
$$ \Omega=S_1Y_1\otimes Q_1+\ldots+S_kY_k\otimes Q_k$$
where $S_i$ are as described just above, $Y_i$ in the center of
$B_i$, and $Q_i\in D$.

An important remark is that we cannot have in our decomposition
both $B_i$ simple and non-simple. Indeed, let us assume, without
loss of generality, that $B_1$ is simple and $B_2$ is not. Let us
recall that the component $B_i$ arise in connection with the
defining tuple of $C$. In this case we consider the first three
components of the tuple $g_1^{(r_1)}$, $(g'_2)^{(r_2)}$ and
$(g''_2)^{(r_2)}$ where all three elements $g_1,g'_2$ and $g''_2$
are pairwise different. Since $B_2$ is not simple, we must have
$e_2=e'_2+e''_2$ where the summands on the right hand side are the
identity elements of the simple components of $B_2$. If we choose
$U=e_1Ue'_2$ then its degree is $g^{-1}_1 g'_2$. If we apply
$\sigma$ to $U$ this element will be mapped into
$e_1S_1Y^{-1}_1(\sigma(U))e''_2Q^{-1}_1Q_2$. If we decompose
$Q^{-1}_1Q_2$ as the sum of homogeneous components of degrees
$t_1,\ldots,t_s$ we obtain that $\sigma(U)$ is the sum of
homogeneous elements of degrees $g^{-1}_1g''_2t_j$. It follows
then that $(g''_2)^{-1}g'_2=t_i$. However, the left hand side of
the latter equation is an element of $\mbox{Supp}\,C$. Therefore,
$g''_2=g'_2$, which is a contradiction.

Now suppose that all $B_i$ in the decomposition of $R_e$ are
simple. Any element of $C$ is the sum of elements of the form
$e_iUe_j$. If we apply $\sigma$ to this element and argue in the
same way as in the previous paragraph then we determine that it is
mapped into $C$ itself. In this case $\sigma(C)=C$, and
$\sigma(D)=D$. The lemma is proved.

\end{proof}

\begin{theorem}
Let $R=M_{n,m}(F)=\oplus_{g\in G} R_g$ be a non-trivial matrix
superalgebra over an algebraically closed field of characteristic
zero or coprime to the order of $G$. Suppose that $*:R\mapsto R$
is a graded superinvolution. Then $R$ as a $G$-graded superalgebra
with superinvolution is isomorphic to the tensor product $C\otimes
D$ where
\par\medskip
\emph{(1)} $C, $ $D$ are graded subsuperalgebras stable under $*$,
and $D$ is a trivial superalgebra;

\emph{(2)} The $G$-grading on $C$ is as in Theorem 5.2 if $*$ is
orthosymplectic on $C$ or as in Theorem 5.3 if $*$ is transpose on
$C$;

\emph{(3)} The $G$-grading on $D$ is as in Theorem 4.1.

\end{theorem}
\begin{proof}

First, we are given that $R=\oplus_{g\in G} R_g$ such that
$(R_g)^*=R_g$ and each $R_g$ is a ${\Bbb Z}_2$-graded vector
space. As before, $R=C\otimes D$ where the $G$-grading on $C$ is
elementary and the $G$-grading on $D$ is fine. Let us consider the
unit component $R_e=C_e\otimes I$ where $I$ is the identity of
$D$. Let $\sigma$ denote an automorphism of $R$ of order 2 that
defines a superalgebra structure on $R$ in the standard way, i.e.
$R_{\bar 0}=\{x\in R|\, \sigma(x)=x\}$ and $R_{\bar 1}=\{x\in R|\,
\sigma(x)=-x\}$. Then $R_e$ is $*$- and $\sigma$-stable. According
to Lemma 6.6, $*$ and $\sigma$ induce respectively a
superinvolution and a superalgebra structure on $C_e$ compatible
with each other.

By Lemma 6.6, $C_e=S_1\oplus\ldots\oplus S_k$ where each $S_i$ is
a $*$-simple superalgebra.  According to Lemma 6.7, either all
$S_i$ in the latter decomposition are non-simple with respect to
$\sigma$, or both $C\otimes I$ and $I\otimes D$ are
$\sigma$-stable.

 If the first possibility holds, then  $C_e$ has a non-trivial ${\Bbb
 Z}_2$-grading.
 Hence we can
choose a non-trivial non-simple $*$-simple $S_i$.  By Lemma 6.2 we
can choose an idempotent $f\in S_i$ such that $f^*=f$ and
$f=f_0+f_1$, $f_1\ne 0$. Set $R'=(f\otimes I)R(f\otimes
I)=fCf\otimes D$. Our main goal is to prove that $fCf\otimes I$ is
$*$- and $\sigma$-stable.
 Note that $fCf\otimes I=(f\otimes I)(C\otimes
I)(f\otimes I)$. Hence, by Lemma 6.3 $(fCf\otimes I)^* \subseteq
fCf\otimes D$. Similarly, $\sigma(fCf\otimes I)\subseteq
fCf\otimes D$.  Let $S_i=A\oplus A^*$ where $A=M_{r,s}(F)$ or
$A=Q(n)$. Denote the identity of $S_i$ by $e_i $. Since
$fCf\subseteq e_iCe_i$, applying Lemma 6.6 (4), we obtain
$(fCf\otimes I)^*\subseteq (e_iCe_i\otimes I)^*=e_iCe_i\otimes I$,
and $\sigma(fCf\otimes I)\subseteq \sigma(e_iCe_i\otimes
I)=e_iCe_i\otimes I$. Therefore, $(fCf\otimes I)^*\subseteq
(e_iCe_i\otimes I)\cap (fCf\otimes D)=fCf\otimes I$. Similarly,
$\sigma(fCf\otimes I)\subseteq (e_iCe_i\otimes I)\cap (fCf\otimes
D)=fCf\otimes I$. Hence, $fCf\otimes I$ is $*$- and
$\sigma$-stable. Note that $Z_{R'}(fCf\otimes I)=f\otimes D\cong
I\otimes D$ as a $G$-graded superalgebras. Therefore,  a
non-trivial superalgebra $f\otimes D$ has a fine grading respected
by a superinvolution, which is a contradiction (see Theorem 4.3).

Next we consider the remaining case when both $C\otimes I$ and
$I\otimes D$ are $\sigma$-stable. If $I\otimes D$ is a non-trivial
subsuperalgebra, then for some $i$, $e_i\otimes D$ is also a
non-trivial superalgebra because $I\otimes D\subseteq e_1\otimes
D+\ldots+e_n\otimes D$.  Hence,  $e_i\otimes D$ is a non-trivial
superalgebra with a fine involution grading, which is wrong (see
Theorem 4.3). The only possible case is when $I\otimes D\subseteq
R_{\bar 0}$. Since $C=C_{\bar 0}\oplus C_{\bar 1}$, $C_{\bar
1}\ne\{0\}$, we have that $R_{\bar 0}=C_{\bar 0}\otimes D$. We
have that $C_{\bar 0}=I_1\oplus I_2$, the sum of two orthogonal
ideals. Therefore, $R_{\bar 0}=(I_1\otimes D)\oplus (I_2\otimes
D)$, the sum of two orthogonal ideals.  It is well known that
 $*$ is either orthosymplectic or  transpose. Directly
applying $*$ to any $I\otimes d$, $d\in D$, we can show that
$I\otimes D$ is $*$-stable. Therefore, $C\otimes I$ is also
$*$-stable. Now 2), 3) follow from Theorems 5.2, 5.3 and 4.1. The
proof is complete.
\end{proof}

\section{Gradings on Lie and Jordan superalgebras}
In this section we apply the results from the previous section in
order to describe group gradings on Jordan superalgebras of types
$osp(n,m)$, $P(n)$ and Lie superalgebras of types $B(n,m)$. As
usual, $H(R,*)$ denotes the set of all elements of a superalgebra
$R$ fixed by a superinvolution $*$ while $K(R,*)$ denotes the set
of all elements of $R$ skew-symmetric under $*$.

The following analog of Proposition 1 (see \cite{BSZ}) holds for
superalgebras. Since the proof of this lemma is a repetition of
that of the just mentioned proposition, we omit the proof.
\begin{lemma}
If $R$ is a $G$-graded superalgebra, then the spaces $H(R,*)$ of
symmetric and $K(R,*)$ of skew-symmetric elements are $G$-graded
if and only if $*$ is a graded involution. $\Box$
\end{lemma}

\begin{theorem}
Let $J$ be a simple Jordan superalgebra of type $osp(n,2m)$ over
an algebraically closed field $F$ of characteristic zero graded by
a finite abelian group $G$. Then there exists a decomposition
$n+2m=2^k q$ such that $J$ as a $G$-graded superalgebra is
isomorphic to the subsuperalgebra $K(R,*)\otimes
K(T,*)+H(R,*)\otimes H(T,*)$ of $R\otimes T$ where $R=M_{2^k}(F)$
is a matrix algebra with an involution  fine grading described in
Theorem 4.1, and $T=M_{s,r}(F) $, $s+r=q$ is a matrix superalgebra
with an involution elementary $G$-grading. Moreover,
 $K(T,*)$ and $H(T,*)$ are as in Theorem 5.2.
\end{theorem}
\begin{proof}
A Jordan superalgebra $J=osp(n,2m)$ is isomorphic to $H(A,osp)$
where $A=M_{n,2m}(F)$ is the matrix superalgebra with the
orthosymplectic superinvolution $osp$. Now let $J=\oplus_{g\in
G}\,J_g$ be a $G$-grading on $J$. Since all automorphisms of
$osp(n,2m)$ can be extended to automorphisms of $M_{n,2m}(F)$, any
grading on $J$ is induced from $A$. It is also easy to observe
that any abelian group of automorphisms is induced by an abelian
group of automorphisms of $A$ because this latter algebra is
generated by its symmetric elements. Also, if $\widehat G$, the
dual group for $G$, leaves invariant the set $H(A,*)$, it also
leaves invariant the set $K(A,*)=[H(A,*),H(A,*)]$. It follows then
by Lemma 7.1 that there exists an involution $G$-grading on
$A=\oplus_{g\in G}\,A_g$ such that $J=H(A,*)$ and $J_g=J\cap A_g$
for all $g\in G$. Clearly, any involution $G$-grading on $A$
induces a $G$-grading on $H(A,*)$.

Now direct computations show that the subspace $H(A
,*)=K(R,*)\otimes K(T,*)+H(R,*)\otimes H(T,*)$ of $R\otimes T$
where $R=M_{2^k}(F)$ is a matrix algebra with an involution fine
grading described in Theorem 4.1, and $T=M_{s,r}(F) $, $s+r=q$ is
a matrix superalgebra with an involution elementary $G$-grading
described in Theorem 5.2. The proof is complete.
\end{proof}

All automorphisms of $P(n)$ are also induced from $M_{n,n}(F)$
which is a special enveloping algebra of $P(n)$. Using the same
arguments we can obtain the description of all gradings on $P(n)$.

\begin{theorem}
Let $J$ be a simple Jordan superalgebra of type $P(n)$ over an
algebraically closed field $F$ of characteristic zero graded by a
finite abelian group $G$. Then there exists a decomposition
$2n=2^k q$ such that $J$ as a $G$-graded superalgebra is
isomorphic to the subsuperalgebra $K(R,*)\otimes
K(T,*)+H(R,*)\otimes H(T,*)$ of $R\otimes T$ where $R=M_{2^k}(F)$
is a matrix algebra with an involution fine grading described in
Theorem 4.1, and $T=M_{s,r}(F) $, $s+r=q$ is a matrix superalgebra
with an involution elementary $G$-grading. Moreover,
 $K(T,*)$ and $H(T,*)$ are as in Theorem 5.3. $\Box$
\end{theorem}
\par\medskip
\par\medskip
According to Serganova's results \cite{Ser}, the only
finite-dimensional simple classical Lie superalgebra that has no
outer automorphisms and can be realized as the set of
skew-symmetric elements is $osp(2n+1,2m)$. Hence, all
automorphisms of $osp(2n+1,2m)$ are induced from $M_{2n+1,2m}(F)$.
The following is an analogue of Theorems 7.2 and 7.3.
\begin{theorem}
Let $L$ be a simple Lie superalgebra of type $osp(2n+1,2m)$ over
an algebraically closed field $F$ of characteristic zero graded by
a finite abelian group $G$. Then there exists a decomposition
$2n+1+2m=2^k q$ such that $L$ as a $G$-graded superalgebra is
isomorphic to the subsuperalgebra $K(R,*)\otimes
H(T,*)+H(R,*)\otimes K(T,*)$ of $R\otimes T$ where $R=M_{2^k}(F)$
is a matrix algebra with an involution fine grading described in
Theorem 4.1, and $T=M_{s,r}(F) $, $s+r=q$ is a matrix superalgebra
with an involution elementary $G$-grading. Moreover,
 $K(T,*)$ and $H(T,*)$ are as in Theorem 5.2.
\end{theorem}

\end{document}